\documentclass[10pt,a4paper,twoside]{article}

\usepackage[utf8,english,pgf,nopscyr]{scinotes}

\UDK{519.22}

\authors
{O.~Volkov, Yu.~Volkov, N.~Voinalovych}
{Волков~О., Волков~Ю., Войналович~Н.}

\author
{Волков Олександр Юрійович}
{студент 1-го року навчання магістерської програми за спеціальністю <<Cтатистика>> Університету Каліфорнії, Берклі}
{статистика, стохастичне моделювання, теорія операторів, машинне навчання}
\author
{Волков Юрій Іванович}
{доктор~фіз.-мат. наук, професор}
{теорія ймовірностей, статистика, дискретна математика}
\author
{Войналович Наталія Михайлівна}
{канд.~пед. наук, доцент кафедри математики, фізики та методик навчання Центральноукраїнського державного університету імені Володимира Винниченка}
{методика навчання математики, теорія ймовірностей, статистика}
\author
{Oleksandr Volkov}
{University of California at Berkeley, Master’s student at the Department of Statistics}
{statistics, stochastic modeling, operator theory, machine learning}
\author
{Yurii Volkov}
{Doctor of Physics an dMathematics, Professor}
{probability theory, statistics, discrete mathematics}
\author
{Nataliia Voinalovych}
{Candidate of Pedagogical Sciences, Associate Professor of the Department of Mathematics and Digital Technologies, Volodymyr Vynnychenko Central Ukrainian State University}
{methodology of teaching mathematics, probability theory, statistics}

\title
[Study of Power Series Distributions with Specified Covariances]
{Study of Power Series Distributions with Specified Covariances}
{Дослідження розподілів степеневих рядів із наперед заданими коваріаціями}

\abstract
{This paper presents a study of power series distributions (PSD) with prescribed covariance characteristics. Such distributions constitute a fundamental class in probability theory and mathematical statistics, as they generalize a wide range of well-known discrete distributions and enable the description of various stochastic phenomena with a predetermined variance structure. 

The aim of the research is to develop analytical methods for constructing power series distributions with given covariances and to establish the conditions under which a particular function can serve as the covariance of a certain PSD. The paper derives a first-order differential equation for the generating function of the distribution, which determines the relationship between its parameters and the form of the covariance function. It is shown that the choice of an analytical or polynomial covariance completely specifies the structure of the corresponding generating function. Based on this, classes of distributions are established that realize linear, quadratic, and more general analytical dependencies between the expectation and the variance.

The analysis made it possible to construct new families of PSDs that generalize the classical Bernoulli, Poisson, geometric, and other distributions while preserving a given covariance structure. The proposed approach is based on the analytical relationship between the generating function and the covariance function, providing a framework for constructing stochastic models with predefined dispersion properties.

The results obtained expand the theoretical framework for describing discrete distributions and open up opportunities for practical applications in statistical estimation, modeling of complex systems, financial processes, machine learning where it is crucial to control the dependence between the mean and the variation. Further research may focus on constructing continuous analogues of such distributions, studying their limiting properties, and applying them to problems of regression and Bayesian analysis.
}
{У роботі проведено дослідження розподілів степеневих рядів (РСР), для яких наперед задано коваріаційні характеристики. Такі розподіли становлять фундаментальний клас у теорії ймовірностей та математичній статистиці, оскільки вони узагальнюють велику кількість відомих дискретних розподілів і дозволяють описати широкий спектр стохастичних явищ з наперед визначеною структурою варіації. 

Метою дослідження є розробка аналітичних методів побудови розподілів степеневих рядів із заданими коваріаціями та встановлення умов, за яких певна функція може виступати коваріацією деякого РСР.  У статті отримано рівняння першого порядку для генератриси розподілу, що визначає взаємозв’язок між її параметрами та формою коваріаційної функції. Показано, що вибір аналітичної або многочленної коваріації повністю визначає структуру відповідної генератриси. На основі цього встановлено класи розподілів, які реалізують лінійні, квадратичні та більш загальні аналітичні залежності між математичним сподіванням і дисперсією.

Проведений аналіз дозволив побудувати нові сімейства РСР, які узагальнюють класичні розподіли Бернуллі, Пуассона, геометричний та інші, зберігаючи при цьому задану форму коваріації. Запропонований підхід базується на аналітичному зв’язку між генератрисою розподілу та функцією коваріації і забезпечує можливість конструювання стохастичних моделей із наперед визначеними дисперсійними властивостями.

Отримані результати розширюють теоретичний апарат опису дискретних розподілів і відкривають можливості для практичного застосування у статистичному оцінюванні, моделюванні складних систем, фінансових процесів, машинному навчанні, де важливо контролювати залежність між середнім і варіацією. Подальші дослідження можуть бути спрямовані на побудову неперервних аналогів таких розподілів, вивчення їхніх граничних властивостей і застосування у задачах регресійного та байєсівського аналізу.}

{розподіли степеневих рядів, генератриса розподілу, коваріаційна функція, дисперсійна характеристика, cтохастична модель}

\begin{document}

\maketitle

\begin{multicols}{2}


\section{ Introduction}
Power series distributions (PSDs) were first systematically considered by Noack~\cite{noack1950} in the context of generalizing classical discrete distributions. Since then, they have taken an important place in probability theory and stochastic modeling, particularly due to their ability to describe a wide class of random variables with specified generating function properties. Further development of the PSD theory is associated with the works of Johnson, Kotz, and Balakrishnan~\cite{johnson1997}, who presented a generalized theory of discrete distributions using power series in their fundamental work. Significant contributions to the study of the properties of such distributions were also made by Morgan~\cite{morgan1984} and Kemp~\cite{kemp1968}, who considered the classification and parameterization of PSDs based on their moment characteristics.

Of particular interest are distributions for which the covariance characteristic (variance function) has a pre-specified form. Similar approaches have found application in statistical modeling where it is necessary to construct stochastic structures with predefined dependencies between the mean and the variance. In this context, the study of distributions with specified covariances is a natural continuation and generalization of the classical results obtained for the Poisson, Bernoulli, and geometric distributions~\cite{ord1967},~\cite{patil1962}.

In recent decades, there has been a growing interest in analytical methods for describing such distributions through equations for their generating functions, which allows for obtaining new classes of probability models with an analytical dependence between parameters~\cite{steutel2003}. In particular, Wolfowitz~\cite{wolfowitz1942} developed a framework for analyzing distributions with specified moment characteristics, which is now also applied to PSDs.

However, the question of constructing PSDs with a specified analytical form of the covariance remains open and is of both theoretical and applied interest. In particular, it is important to determine the conditions under which an arbitrary analytical or polynomial function can act as the covariance of a certain power series distribution.

The relevance of the research is due to the necessity of expanding the class of known PSDs by introducing a parameterization based on specified covariances. This approach allows for constructing new types of probability distributions with controlled variance properties, which can be useful in problems of stochastic modeling, statistical estimation, and machine learning. Furthermore, establishing a connection between the functional form of the covariance and the corresponding generating function opens possibilities for further classification of PSDs and generalization of known results.

The goal of the work is to investigate the properties of power series distributions with specified covariances, establish the equation satisfied by the generating function of such a distribution, and determine the conditions under which an arbitrary function can be the covariance of a PSD.

\section{Power Series Distributions with Specified Covariances}


The term  ``power series distributions'' (PSD) belongs to Noack~\cite{noack1950}. These distributions are defined as follows: let a power series with non-negative coefficients be given
$$
\omega(y) = \sum_{k=0}^{\infty} a_k y^k, ~a_k \ge 0, ~0 \le y < R,
$$
($R$ is the radius of convergence of the series). Consider a random variable $\xi$ that can take non-negative integer values with probabilities
$$
p_k = P\{\xi = k\} := a_k \frac{y^k}{\omega(y)}, ~k=0, 1, 2, \ldots
$$
Such a distribution is called the power series distribution of the function $\omega(y)$. Its probability generating function is
$$
P(z) = \sum_{k=0}^{\infty} a_k \frac{y^k}{\omega(y)} z^k = \frac{\omega(yz)}{\omega(y)}.
$$

From this, we obtain the expected value and the variance of $\xi$.
We have
$$
\text{M}\xi = P'(1) = y \frac{\omega'(y)}{\omega(y)},
$$
$$
\text{D}\xi = P''(1) + P'(1) - (P'(1))^2,
$$
and hence
$$
\text{D}\xi = y^2 \frac{\omega''(y)}{\omega(y)} + y \frac{\omega'(y)}{\omega(y)} - y^2 \left(\frac{\omega'(y)}{\omega(y)}\right)^2.
$$
Let's introduce a new parameterization

\begin{equation}\label{eq0}
x = y \frac{\omega'(y)}{\omega(y)}.
\end{equation}
Then
$$
\text{D}\xi = y \frac{dx}{dy}.
$$

The function $x=x(y)$, which is defined by the relation~\eqref{eq0}, has an inverse on the interval $(0,R)$, because $\frac{dx}{dy} = \frac{\text{D}\xi}{y} > 0$. If we denote this inverse function by $y=f(x)$, then with the new parameterization, the variance $\text{D}\xi$ will be a function of $x$. This function we will denote by $V(x)$ and call the covariance characteristic (covariance) of the distribution $\xi$. Thus,
$$
V(x) = \frac{f(x)}{f'(x)}.
$$
After introducing this parameterization, the generating function will be a function of two variables $P=P(z,x)$.

\vspace{0.1cm}
\begin{theorem} \label{th1}
The generating function satisfies the following first-order equation:
\begin{equation}\label{eq1}
V(x) \frac{\partial P}{\partial x} - z \frac{\partial P}{\partial z} + xP = 0,
\end{equation}
$P(1,x)=1,$ $\forall~x \in X$, where $X$ is the range of the function $x=x(y)$.
\end{theorem}

\vspace{0.1cm}

\begin{proof}
Let the probability generating function of the power series distribution be given as $P(z) = \frac{\omega(yz)}{\omega(y)}$.
We consider the definition of the parameterization $x$ and the covariance characteristic $V(x)$:
$$
x = y \frac{\omega'(y)}{\omega(y)},~\text{where } y=f(x);~V(x) = \frac{f(x)}{f'(x)}.
$$
Then we obtain the relationship between the differentials:
$$
f'(x) = \frac{dy}{dx} \implies V(x) = \frac{y}{dy/dx} \implies V(x) \frac{dy}{dx} = y.
$$
We apply the chain rule for $\frac{\partial P}{\partial x}$:
$$
V(x) \frac{\partial P}{\partial x} = V(x) \frac{\partial P}{\partial y} \frac{dy}{dx} = y \frac{\partial P}{\partial y}.
$$
We compute the partial derivatives $\frac{\partial P}{\partial y}$ and $\frac{\partial P}{\partial z}$:
$$
\frac{\partial P}{\partial y} = \frac{\partial}{\partial y} \left( \frac{\omega(yz)}{\omega(y)} \right) = \frac{z\omega'(yz)\omega(y) - \omega'(y)\omega(yz)}{\omega^2(y)},
$$
$$
\frac{\partial P}{\partial z} = \frac{\partial}{\partial z} \left( \frac{\omega(yz)}{\omega(y)} \right) = \frac{y\omega'(yz)}{\omega(y)}.
$$

We substitute the obtained expressions into the left side of equation~\eqref{eq1}:

$$
V(x) \frac{\partial P}{\partial x} - z \frac{\partial P}{\partial z} + xP = y \frac{\partial P}{\partial y} - z \frac{\partial P}{\partial z} + xP =
$$

$$
= y \left(\frac{z\omega'(yz)\omega(y) - \omega'(y)\omega(yz)}{\omega^2(y)}\right) -z \left(\frac{y\omega'(yz)}{\omega(y)}\right) +
$$
$$
+ \left(y \frac{\omega'(y)}{\omega(y)}\right)\cdot \left(\frac{\omega(yz)}{\omega(y)}\right)=\frac{yz\omega'(yz)\omega(y)}{\omega^2(y)} -
$$
$$
-\frac{y\omega'(y)\omega(yz)}{\omega^2(y)}- \frac{yz\omega'(yz)}{\omega(y)} + \frac{y\omega'(y)\omega(yz)}{\omega^2(y)}=
$$
$$
= \frac{yz\omega'(yz)}{\omega(y)} - \frac{y\omega'(y)\omega(yz)}{\omega^2(y)} - \frac{yz\omega'(yz)}{\omega(y)} +
$$
$$
+\frac{y\omega'(y)\omega(yz)}{\omega^2(y)} = 0.
$$
The equation is satisfied.
\end{proof}

\vspace{0.1cm} %

\begin{example}
$\omega(y)=1+y$.
\\
In this case, we obtain the Bernoulli distribution with parameter $p=x$,
$$
x=\frac{y}{1+y},~y=\frac{x}{1-x}, ~P(z,x)=1+x(z-1),
$$
$$
V(x)=x(1-x), ~0<x<1
$$
\end{example}

\vspace{0.1cm} %

\begin{example}
$\omega(y)=\frac{1}{1-y}$.
\\
In this case, we obtain the geometric distribution with parameter $p=y$
$$
x=\frac{y}{1-y}, ~y=\frac{x}{1+x}, ~P(z,x)=\frac{1}{1+x(1-z)},
$$
$$
V(x)=x(1+x), ~x>0.
$$
\end{example}

\vspace{0.1cm} %

\begin{example}
$\omega(y)=\exp y$.
\\
In this case, we obtain the Poisson distribution with parameter $\lambda=x,~x=y,$
$$
P(z,x)=\exp(x(z-1)),~V(x)=x.
$$
\end{example}

We note that in the given examples, the covariance is either a linear function of $x$ or a quadratic one. As shown in~\cite{morris1982}, other PSDs for which the covariance is a polynomial of degree no higher than 2 do not exist (up to transformations, which will be discussed below). In the work~\cite{volkov1990}, all distributions for which the covariance is a cubic polynomial were found.

A natural question arises: how do we search for PSDs for which the covariance will be an a priori specified function, for example, a polynomial? To answer this question, we establish several auxiliary statements.

\vspace{0.1cm} %
\noindent\textbf{Lemma 1.}
\textit{
Let the function $V(x)$ be the covariance of the PSD of the function $\omega(y)$. Then
for all $m \in \mathbb{N}$ the function 
$V(x-m)$ will be the covariance of the PSD of the function
$\omega_m(y)=C y^m \omega(y)$, $C$--is an arbitrary positive constant.
}
\vspace{0.1cm} %
\begin{proof}
We introduce the parameter $x$ by the formula
$$
x=y \frac{\omega'_m(y)}{\omega_m(y)}.
$$
Hence
$$
x = m+y \frac{\omega'(y)}{\omega(y)},
$$
therefore $ y=f(x-m),$ thus
$$
\frac{y(x)}{y'(x)}=V(x-m).
$$
\end{proof}

Analogously, it is proved

\vspace{0.1cm} %
\noindent\textbf{Lemma 2.}\textit{~Let the function $V(x)$ be the \text{covariance} of the PSD of the function $\omega(y)$ and
$$
\omega(0)=\omega'(0)=\dots=\omega^{(r-1)}(0), ~\text{and } ~\omega^{(r)}(0) \ne 0.
$$
Then the function $V(x+r)$ will be the covariance of the PSD of the function $C y^{-r} \omega(y)$, $C$ is an arbitrary positive constant.
}

\vspace{0.1cm} %

\noindent\textbf{Lemma 3.}\textit{~Let the function $V(x)$ be the covariance of the PSD of the function $\omega(y)$. Then
$\forall n \in \mathbb{N}$ the function $nV(x/n)$ will be the covariance of the PSD of the function $C(\omega(y))^n$,
$C$ is an arbitrary positive constant.
}
\vspace{0.1cm}
\begin{proof}
We introduce the parameter $x$ by the formula
$$x=y\frac{d}{dy}\left(\log(\omega(y))^n\right),$$
hence
$$
x = n y \frac{\omega'(y)}{\omega(y)},
$$
therefore $y=f(x/n)$, thus,
$$ \frac{y(x)}{y'(x)}=nV\left(\frac{x}{n}\right).$$
\end{proof}

\vspace{0.1cm} 
\noindent\textbf{Lemma 4.}\textit{~Let the function $V(x)$ be the covariance of the PSD of the function $\omega(y)$. Then $\forall n\in\mathbb{N}$ the function $n^{2}V\left(\frac{x}{n}\right)$ will be the covariance of the PSD of the function $C\omega(y^n)$, $C$ --- is an arbitrary positive constant.
}
\vspace{0.1cm}
\begin{proof}
We introduce the parameter $x$ by the formula
$$
x = y\frac{d}{dy}\left(\log\omega(y^n)\right),~\text{hence} ~x = ny^{n}\frac{\omega'(y^n)}{\omega(y^n)}
$$
therefore $y^n=f(x/n)$,
thus,
$$\frac{y(x)}{y'(x)}=n^{2}V\left(\frac{x}{n}\right)$$.
\end{proof}

\vspace{0.1cm}

\noindent\textbf{Corollary 1.} \textit{Let the function $V(x)$ be the covariance of the PSD of the function $\omega(y)$. Then
$\forall m \in \mathbb{Z}, n \in \mathbb{N}, k \in \mathbb{N}$ the function
$$
k n^2 V\left(\frac{x-m}{k n}\right)
$$
will be the covariance of the PSD of the function $C y^m (\omega(y^n))^k$, $C$ – is an arbitrary positive constant.
}

\vspace{0.1cm}

Arbitrary shifts and scaling of the covariance will not, generally speaking, be the covariance of a PSD; for example, the function $2x(1-x)$ cannot be a covariance. Thus, another problem arises: what properties must the function $V(x)$ possess for it to be the covariance of a PSD. This question can be answered as follows. Let $V(x)$ be an arbitrary continuous function on some set $X$. Consider equation~\eqref{eq1}. Such an equation with initial conditions $P(1,x)=1,~\forall x\in X$ has a unique solution, which can be found in quadratures. If the obtained solution is the generating function of some PSD, then the function $V(x)$ will be the covariance.

Let's find the solution to equation~\eqref{eq1}. For this, we introduce the following notation:
$$
s(z) := \exp\left(\int \frac{dz}{V(z)}\right), ~b(z) := \exp\left(\int \frac{z dz}{V(z)}\right),
$$
where
$$\int \frac{dz}{V(z)} \quad\text{and} \quad\int \frac{z dz}{V(z)}$$ are arbitrary antiderivatives of the respective functions.

Then $P(z,x)=(b(x))^{-1}b(s^{-1}(z \cdot s(x)))$, where $s^{-1}$ is the inverse function of $s$ (if such a function exists). Next, we need to investigate the function $P(z,x)$. A necessary condition for this function to be a generating function is its analyticity in the neighborhood of the point $z=0$. This condition holds when the function $V(z)$ is analytic in some circle $\{z:|z|<R\}$, takes nonzero values on the interval $(0,R)$ and $\tau(0) \ne 0$, where $\tau(z)=\frac{z}{s(z)}$. By Lagrange's theorem~\cite[p.~80]{whittaker1963}
$$
P(z,x) = (b(x))^{-1} \sum_{k=0}^{\infty} c_k (s(x))^k z^k, ~\text{where} ~c_0=b(0),
$$
\begin{equation}\label{eq3}
c_k = \frac{1}{(k+1)!}\left(\frac{d}{dz}\right)^k (b'(\tau(z))(\tau(z))^{k+1})\Big|_{z=0},
\end{equation}
$ k=1, 2, \dots$

Therefore, if $s(x)>0, c_k \ge 0$, then the function $P(z,x)$ will be the generating function, and $V(x)$ – the covariance of the PSD of the function $\omega(y)=Cb(s^{-1}(y))$, $C$ – is an arbitrary positive constant.

Thus, the following is proved:

\vspace{0.1cm} %
\noindent\textbf{Theorem 2.}\textit{~If the function $V(z)$ is analytic in some circle $\{z:|z|<R\}$, takes positive values on the interval $(0,R)$, $\tau(0)>0$ and the coefficients $c_k$, found by formulas~\eqref{eq3}, are non-negative, then the function $V(x)$ is the covariance of the PSD of the function $\omega(y)=b(z)$, where $y=s(z)$.
}
\vspace{0.1cm} %

\vspace{0.1cm} %

\begin{example}
$$
V(x)=x(1+x(\alpha-1))(1+x\alpha), ~\alpha>1, ~x>0;
$$
$$
b(z)=\frac{1+z\alpha}{1+(\alpha-1)z}, ~s(z)=z \frac{(1+(\alpha-1)z)^{\alpha-1}}{(1+z\alpha)^\alpha}.
$$
If we replace $z$ with $\frac{z}{1-(\alpha-1)z}$, then
$$
\omega(y)=\sum_{k=0}^{\infty} c_k y^k=1+z, ~\text{where}
$$
$$
y=\frac{z}{(1+z)^\alpha},~0 \le z < 1, ~c_k = \frac{1}{\alpha k+1} \binom{\alpha k+1}{k}.
$$
\end{example}
When using Theorem 2 to find a PSD with a given covariance, the most difficult part is checking the signs of the coefficients $c_k$.

\vspace{0.1cm} %
\noindent\textbf{Theorem 3.}\textit{~If $V(x)=x(1+xU(x))$, $U(x)$ is absolutely monotonic on
some interval $(0,R)$, then the function $V(x)$ is the covariance of the PSD of the function $b(s^{-1}(y))$.
}

\vspace{0.1cm} %

\begin{proof}
Let
$$U(x)=a_1+a_2 x+\dots. $$
We introduce an auxiliary function
$$
g(z)=b'(z)(\tau(z))^{a+1}=\frac{z b(z)}{V(z)}(\tau(z))^{a+1}=
$$
$$
=\sum_{k=0}^{\infty} t_k z^k, \quad t_k=t_k(a).
$$
Hence
\\
$z V(z) g'(z) =$
$$ = ((a+2)V(z) - z V'(z) + z^2 - (a+1)z)g(z),$$
and since
$$g'(z)=\sum_{k=0}^{\infty} k t_k z^{k-1},$$
for the coefficients $t_k$ we obtain the following recurrence relation
\\
$(k+1)t_k=(a_1(a-k)+1)t_{k}+(a-k)(a_2 t_{k-1}+a_3 t_{k-2}+\dots+a_{k+1} t_0)$, $t_{-i}=0$ for $i>0$, $t_0=1$. If $a=m$, we sequentially find $t_1(m), t_2(m), \dots, t_m(m)$, and then $$c_m = \frac{t_m(m)}{m+1}, ~c_0=1.$$
These coefficients are non-negative and will be the coefficients of the power series of the function $b(s^{-1}(y))$, which was required to be proved.
\end{proof}

\vspace{0.1cm} %

\noindent\textbf{Corollary 2.}\textit{~If $U(x)$ --- is a polynomial with non-negative coefficients, then the function $V(x)=x(1+x U(x))$ will be the covariance of a PSD.
}

\vspace{0.1cm} %

Other polynomial covariances can be obtained from $V(x)$ by using Corollary 1.

\vspace{0.1cm} %

\begin{example}
$$
V(x)=x(1+2x)(1+3x)(1+4x), \ x>0,
$$
$$
b(z)=\frac{(3z+1)^{3}(6z+1)^{3/2}}{(4z+1)^{4}\sqrt{2z+1}},
$$
$$
s(z)=\frac{z(2z+1)(4z+1)^{16}}{(3z+1)^{9}(6z+1)^{9}},
$$
$
c_{0}=c_{1}=1, ~c_{2}=8, ~c_{3}=96, ~c_{4}=1379, ~c_{5}=21937, ~c_{6}=372724, ~c_{7}=6631164, ~\dots
$
\end{example}

\vspace{0.1cm} %

\begin{example}
$$V(x)=x\left(1+x/2\right) \left(1+x\right)^{2}; \ x>0,
$$
$$
b(z)=\left(\frac{2+z}{2(1+z)}\right)^{2} \exp\left(\frac{2z}{1+z}\right),
$$
$$
s(z)=\frac{z}{2+z}\exp\left(-\frac{2z}{1+z}\right),
$$
$
c_{0}=1, ~c_{1}=2, ~c_{2}=7, ~c_{3}=92/3, ~c_{4}=455/3, ~c_{5}=4046/5, ~c_{6}=204631/45, ~\dots
$
\end{example}

We use Lemma 1 and substitute the expression $x-1$ for $x$ in the formula for $V(x)$. We obtain the covariance $x^2(x^2-1)/2$ of the PSD of the function
$$\omega(y)=\frac{z(z+2)}{4(z+1)^{2}}, ~\text{where} ~y=\frac{z}{2+z}\exp\left(-\frac{2z}{1+z}\right).$$
If we replace $\frac{z}{1+z}$ with $2z$, then $\omega(y)=z(1-z)$, where
$$
y=\frac{z}{1-z}\exp(-4z),
$$
and then the coefficients $c_k$ can also be found by the following formulas: $c_{k+1}=$
$$
\frac{1}{(k+1)!}\left(\frac{d}{dz}\right)^{k}\left[(1-2z)(1-z)\exp(4z)\right]^{k+1}\Big|_{z=0},
$$
$k=0, 1, 2, ~\dots$

\vspace{0.1cm} %

\begin{example}
$
V(x)=x\left(1+x^{3}\right); \ x>0;
$
$$
b(z)=\frac{\sqrt[3]{1+z}}{\sqrt[6]{1-z+z^{2}}}\exp\left(\frac{1}{\sqrt{3}}\operatorname{arctg}\frac{z\sqrt{3}}{2-z}\right),
$$
$$
s(z)=\frac{z}{\sqrt[3]{1+z^{3}}},
$$
$
c_{0}=c_{1}=1, ~c_{2}=1/2, ~c_{3}=1/6, ~c_{4}=1/8, ~c_{5}=11/120, ~c_{6}=31/720, c_{7}=11/240, ~\dots
$
\end{example}

\vspace{0.1cm} %

\begin{example} 
$$
V(x)=x\left(1+x/2\right) \left(1+x+x^{2}/2\right), \ x>0,
$$
$$
b(z)=\frac{1}{2}\frac{(z+2)^{2}}{z^{2}+2z+2}\exp\left(2\operatorname{arctg}\frac{z}{z+2}\right), 
$$
$$
s(z)=\frac{z}{z+2}\exp\left(-2\operatorname{arctg}\frac{z}{z+2}\right),
$$
$
c_{0}=1, ~c_{1}=2, ~c_{2}=5, ~c_{3}=44/3, ~c_{4}=143/3, ~c_{5}=166, ~c_{6}=5465/9, ~\dots
$
\end{example} 

We use Lemma 1 and substitute the expression $x-1$ for $x$ in the formula for $V(x)$. We obtain the covariance $(x^{4}-1)/4$ of the PSD of the function
$$
\omega(y)=\frac{z(z+2)}{2(z^{2}+2z+2)}, ~\text{where}
$$

$$
y=\frac{z}{z+2}\exp\left(-2\operatorname{arctg}\frac{z}{z+2}\right).
$$
If we replace $\frac{z}{z+2}$ with $z$, then $\omega(y)=\frac{z}{1+z^2}$, where
$$
y=\frac{z}{\exp(2\operatorname{arctg}z)},
$$
and then the coefficients $c_k$ can also be found by the following formulas: $c_{k+1}=$
$$
\frac{1}{(k+1)!}\left(\frac{d}{dz}\right)^{k}\left[\frac{1-z^{2}}{(1+z^{2})^{2}}\exp(2\operatorname{arctg}z)\right]^{k+1}\Big|_{z=0},
$$
$k=0, 1, 2, ~\dots$

\vspace{0.1cm} %

\begin{example} 
$
V(x)=x(1+x)^{3}; \ x>0;
$
$$
b(z)=\exp\left(\frac{1}{2}\left(1-\frac{1}{(1+z)^{2}}\right)\right), 
$$
$$
s(z)=\frac{z}{1+z}\exp\left(-\frac{z}{1+z}+\frac{1}{2}\left(\frac{1}{(1+z)^{2}}-1\right)\right),
$$
$
c_{0}=1, ~c_{1}=1, ~c_{2}=2, ~c_{3}=31/6, ~c_{4}=61/4, ~c_{5}=5861/120, ~c_{6}=7438/45, \dots
$
\end{example} 

If we replace $\frac{z}{1+z}$ with $z$, then $$
\omega(y)=\exp\left(z-\frac{z^{2}}{2}\right), ~\text{where} ~y=\exp\left(2z-\frac{1}{2}z^{2}\right),
$$
and then the coefficients $c_k$ can also be found by the following formulas:
\begin{flalign*} c_{k+1}&=\frac{1}{(k+1)!}\times
& \end{flalign*}
$$
\left(\frac{d}{dz}\right)^{k}\left[(1-z)\exp\left((2k+3)z-z^{2}(1+k/2)\right)\right]\Big|_{z=0}
$$
$k=0,1,2, ~\dots$

\vspace{0.1cm} %

\begin{example}
$
V(x)=\frac{x}{1-x}, \ 0<x<1,
$
$$
b(z)=\exp\left(z-\frac{z^{2}}{2}\right), \quad s(z)=z\exp(-z),
$$
$
c_{0}=1, ~c_{1}=1, ~c_{2}=1, ~c_{3}=7/6, ~c_{4}=19/12, ~c_{5}=97/40, ~c_{6}=733/180, ~\dots
$
\end{example}

\vspace{0.1cm} %
\begin{example}
$V(x)=2(1-\sqrt{1-x})$, $0<x<1,$
$$b(z)=\exp\left(\frac{1}{3}\left(1-\sqrt{(1-z)^3}\right)+\frac{1}{2}z\right),$$
$$s(z)=\frac{z}{1+\sqrt{1-z}}\exp(1-\sqrt{1-z}),$$
$c_0=1, ~c_1=2, ~c_2=5/2,~c_3=8/3, ~c_4=67/24, ~c_5=31/10, ~c_6=2731/720, ~\dots$
\end{example}

\section{Conclusion and further developments}
In this work, Power Series Distributions with specified covariances were investigated. A first-order equation for the generating functions of such distributions was obtained, and it was shown that the analytical form of the covariance determines the functional form of the generating function. Conditions were established under which an arbitrary analytical or polynomial function can be the covariance of some Power Series Distribution.

A generalized scheme for finding distributions with given dispersion characteristics was constructed, based on the relationship between the parameterization of the generating function and the covariance function. The proposed approach allowed not only to reproduce classical cases (Bernoulli, Poisson, geometric distributions) but also to obtain new examples of distributions with given polynomial covariances.

The obtained results can be used for constructing stochastic models with controlled variational properties, in statistical machine learning problems, and also when studying generalized discrete distributions in mathematical statistics. Further research can be aimed at analyzing the properties of continuous analogs of such distributions and studying their limiting forms.


\end{multicols}

	


\end{document}